\documentclass{article}
\usepackage{amsfonts}
\title{\large \bf Some Applications of Augmentation Quotients}
\author{\small \bf Deepak Gumber \\
\small \em School of Mathematics and Computer Applications\\
\small \em Thapar University, Patiala - 147 004,
India\\
\small E-mail: dkgumber@yahoo.com}
 
\date{}
\newtheorem{thm}{Theorem}[section]
\newtheorem{lm}[thm]{Lemma}
\newtheorem{pp}[thm]{Proposition}

\begin{document}

\maketitle
\begin{abstract}
We give some applications of augmentation quotients of free group rings in group theory. 
\end{abstract}

\vspace{2ex}


\vspace{2ex}


\section{\large Introduction}

Let $\mathbb{Z}G$ denote the integral group ring of a group $G$ and $\Delta(G)$ its augmentation ideal. Let $\{\gamma_n(G)\}_{n\geq 1}$ be the lower central series of $G$. We also write $G'$ for $\gamma_2(G)=[G,G]$, the derived group of $G$. When $G$ is free, then integral group ring is known as free group ring. Let $\Delta^n(G)$
denote the $n$-th associative power of $\Delta(G)$ with $\Delta^0(G)=\mathbb{Z}G$. The additive abelian group $\Delta^n(G)/\Delta^{n+1}(G)$ is known as the $n$-th augmentation quotient and has been intensively studied during the last forty years. Vermani[7] has given a notable survey article about work done on augmentation quotients. In this short note we are interested in the applications of augmentation quotients in group theory.
Henceforth, unless or otherwise stated, $F$ is a free group and $R$ is a normal subgroup of $F$. Hurley and Sehgal[4] identified the subgroup $F\cap (1+\Delta^2(F)\Delta^n(R))$ for all $n\geq 1$ and then using the fact that $\Delta(F)\Delta^n(R)/\Delta^2(F)\Delta^n(R)$ is free abelian for all $n\geq 1$ [1], they showed that the group $\gamma_{n+1}(R)/\gamma_{n+2}(R)\gamma_{n+1}(R\cap F')$ is a free abelian group for all $n\geq 1$.  Gruenberg [1, Lemma III.5] proved that $\Delta^n(F)\Delta^m(R)/\Delta^{n+1}(F)\Delta^m(R)$ is a free abelian group for all $m,n \geq 1$. When $R$ is an arbitrary subgroup of $F$, Karan and Kumar [5] proved that the groups $\Delta^n(F)\Delta^m(R)/\Delta^{n+1}(F)\Delta^m(R)$, 
$\Delta^n(F)\Delta^m(R)/\Delta^{n-1}(F)\Delta^{m+1}(R)$  and $\Delta^n(F)\Delta^m(R)/\Delta^n(F)\Delta^{m+1}(R)$ are free abelian for all $m,n \geq 1$. They  gave the complete description of all these groups and explicit bases of first two groups. As a consequence of their results they proved that $R'/[R',R\cap F']$ is a free abelian group. Gumber et. al. [2] proved that $\Delta^p(R)\Delta^n(F)\Delta^q(R)/\Delta^p(R)
\Delta^{n+1}(F)\Delta^q(R)$ is free abelian for all $p,q,n\geq 1$ and as a consequence showed that   
$\gamma_3(R)/\gamma_4(R)[R\cap F',R\cap F',R]$ is a free abelian group. In section 3, we identify the subgroup 
$$R\cap (1+\Delta^{m+3}(R)+\Delta^{m+1}(R)\Delta(R\cap F')+\Delta^m(R)\Delta([R,R\cap F']))$$ for $m=0,1,$ and $2$, 
 and then prove\\

\noindent {\bf Theorem A} The groups $R'/[R,R\cap F'],\;
\gamma_3(R)/\gamma_4(R)[R,R\cap F',R\cap F']$, and
$\gamma_4(R)/\gamma_{5}(R)[R,R\cap F',R\cap F',R\cap F'][\,[R,R\cap F'],[R,R\cap F']\,]$
are free abelian.

\section{\large Preliminaries}

Let $G$ be a group and $H$ be a normal subgroup of $G$ such that $G/H$ 
is  free-abelian. Let 
$\{{x_{\delta }H\mid \delta\in \Delta}\}$ be a basis for $G/H$. 
We may suppose that the index set $\Delta $ is well ordered.
As $G^{\prime}\subset H$, $S$,
the set consisting of elements of the form $x_{\delta_{1}}^
{t_{1}}
x_{\delta_{2}}^{t_{2}}\ldots x_{\delta_{n}}^
{t_{n}},\;t_{i}\in Z,n\geq 1,\delta_{1}< \delta_{2}< \ldots < \delta_{n},$
is a transversal of $H$ in $G$. 
Let $L_{n}$  be the $Z$-submodule
of $\Delta(G)$ generated by  elements  of the form
\[(x_{\delta_{1}}^{\epsilon_{1}}-1)\ldots
(x_{\delta_{n}}^{\epsilon_{n}}-1),\;
\epsilon_{i}=1\;\mathrm{or}\;-1\;\mathrm {for\; every}\;i\;\mathrm{and}\;\delta_{1}
\leq \delta_{2}\leq \ldots \leq \delta_{n}.\]
For $m\geq 2$, let $L^{(m)}=\displaystyle\sum_{n\;\geq \;m}L_{n}$.

\begin{thm}
{\em [8]} For $n\geq 2$, $\Delta^{n}(G)$ is equal to
$$\Delta^{n-1}(G)\Delta(H)+
\Delta^{n-2}(G)\Delta(G^{\prime})+\cdots
+\;\Delta(G)\Delta(\gamma_{n-1}(G))+
\Delta(\gamma_{n}(G))\;\oplus\; L^{(n)}.$$
\end{thm}

Let $U$ be a group and $W$ be a left transversal  of a subgroup $V$ 
of  $U$  in $U$ with  $1\in W$. Then every element of $U$ can be uniquely
written as wv,  $w\in W$, $v\in V$. Let  $\phi :ZU\rightarrow ZV$ be the
onto homomorphism of right $ZV$-modules which on the elements of $U$ is
given by $\phi (wv)=v,\;w\in W,\;v\in V$. The homomorphism $\phi$ maps
$\Delta(U)J$ onto $\Delta(V)J$ for every ideal  $J$  of $ZV$. In particular,
by the choice of the transversal $S$  of  $H$
in  $G$, we have $\phi\mid_
{\textstyle L^{(n)}}=0$. The homomorphism $\phi$ is usually called the filtration map.\\

We shall also need the following results:
\begin{lm}{\em [9]}
 Let $G$ be a group, $K$ a subgroup of $G$, and $J$ an ideal of $\mathbb{Z}G$ containing $\Delta^2(K)$. Then
$G\cap (1+J+\Delta(K))=(G\cap (1+J))K.$
\end{lm}
\begin{thm}{\em [8]} Let $G$ be a group with a normal subgroup $H$ such that $G/H$ is free abelian. Then $G\cap (1+\Delta^n(G)+\Delta(G)\Delta(H))=\gamma_n(G)H'$ for all $n\geq 1$.
\end{thm}

\section{\large Proof of Theorem A}
To avoid repeated and prolonged expressions, we write
\begin{eqnarray*}
A&=&\Delta^{4}(R)+\Delta^{2}(R)\Delta(R\cap F')+\Delta(R)
\Delta([R,R\cap F'])\\
B&=&\Delta^{5}(R)+\Delta^{3}(R)\Delta(R\cap F')+\Delta^{2}(R)
\Delta([R,R\cap F']).
\end{eqnarray*}

\begin{pp}
$R\cap (1+\Delta^3(R)+\Delta(R)\Delta(R\cap F')+\Delta([R,R\cap F']))=[R,R\cap F']$.
\end{pp}
{\bf Proof.} Proof is easy and follows by Lemma 2.2 and Theorem 2.3.

\begin{pp} $R\cap (1+A)=\gamma_{4}(R)[R,R\cap F',R\cap F']$.
\end{pp}
{\bf Proof.} Since $\gamma_4(R)-1\subset \Delta^4(R)$ and 
$[R,R\cap F',R\cap F']-1\subset \Delta^{2}(R)\Delta(R\cap F')+\Delta(R)
\Delta([R,R\cap F'])$, it follows that $\gamma_{4}(R)[R,R\cap F',R\cap F']\subset R\cap (1+A)$. For the reverse inequality, we let $w\in R$ such that
$w-1\in A$ and proceed to show that $w\equiv 1\pmod {\gamma_{4}(R)[R,R\cap F',R\cap F']}.$
Since $R/R\cap F'$ is free-abelian, using Theorem 2.1 repeatedly
we have
\begin{eqnarray*}
A&=&\Delta(\gamma_{4}(R))+L^{(4)}+\Delta(R)\Delta^{2}(R\cap F')+
\Delta(R')\Delta(R\cap F')\\
&&\;\;\;\;\;\;\;\;\;\;\;\;\;\;\;\;+\;L^{(2)}\Delta(R\cap F')
+\Delta(R)\Delta([R,R\cap F']).
\end{eqnarray*}
Now since $R\cap(1+A)\subset R\cap F'$, using the filtration map
$\phi : ZR\rightarrow Z(R\cap F')$, it follows that
\begin{eqnarray*}
R\cap(1+A)&\subset&(R\cap F')\cap (1+\Delta^{3}(R\cap F')+\Delta(R')\Delta(R\cap F')\\
&&\;\;\;\;\;\;\;\;\;\;+\Delta(R\cap F')\Delta([R,R\cap F'])
+\Delta(\gamma_{4}(R)))\\
&\subset&(R\cap F')\cap (1+\Delta^{3}(R\cap F')+\Delta(R')\Delta(R\cap F')\\
&&\;\;\;\;\;\;\;\;\;\;+\Delta([R,R\cap F',R\cap F'])+\Delta(\gamma_{4}(R)))\\
&=&(R\cap F')\cap (1+\Delta^{3}(R\cap F')+\Delta(R')\Delta(R\cap F'))\\
&&\;\;\;\;\;\;\;\;\;\;\;\;\;\;\;[R,R\cap F',R\cap F']\gamma_{4}(R)\\
&=&[R,R\cap F',R\cap F']\gamma_{4}(R),
\end{eqnarray*}
where last equality follows by Theorem 2.4 and second last equality follws by Lemma 2.3.
\hfill $\Box$

\begin{pp} $$R\cap (1+B)=\gamma_{5}(R)[R,R\cap F',R\cap F',R\cap F']
[\,[R,R\cap F'],[R,R\cap F']\,].$$
\end{pp}
{\bf Proof.} As in the above proposition, it is sufficient to prove that if $w\in R$ is such that $w-1\in B$, then $$w\equiv 1 \pmod {\gamma_{5}(R)[R,R\cap F',R\cap F',R\cap F']
[\,[R,R\cap F'],[R,R\cap F']\,]}.$$
Using Theorem 2.1 repeatedly, we have
\begin{eqnarray*}
&&\Delta^{5}(R)+\Delta^{3}(R)\Delta(R\cap F')+\Delta^{2}(R)
\Delta([R,R\cap F'])\\
&=&\Delta(R)\Delta(\gamma_{4}(R))+\Delta(\gamma_{5}(R)+
L^{(5)}+\Delta(R)\Delta^{3}(R\cap F')\\
&&\;\;\;\;\;+\Delta(R')\Delta^{2}(R\cap F')
+L^{(2)}\Delta^{2}(R\cap F')+\Delta(R)\Delta(R')\Delta(R\cap F')\\
&&\;\;\;\;\;+\Delta(\gamma_{3}(R))\Delta(R\cap F')
+L^{(3)}\Delta(R\cap F')
+\Delta(R)\Delta(R\cap F')\\
&&\;\;\;\;\;\Delta([R,R\cap F'])
+\Delta(R')\Delta([R,R\cap F'])+L^{(2)}\Delta([R,R\cap F']).
\end{eqnarray*}
Applying filtration map $\phi : ZR\rightarrow Z(R\cap F')$, we have
\begin{eqnarray*}
&&R\cap (1+\Delta^{5}(R)+\Delta^{3}(R)\Delta(R\cap F')+\Delta^{2}(R)
\Delta([R,R\cap F']))\\
&=&(R\cap F')\cap  (1+\Delta^{4}(R\cap F')
+\Delta(R')\Delta^{2}(R\cap F')
+\Delta(\gamma_{3}(R))\Delta(R\cap F')\\
&&\;\;\;\;\;+\Delta^{2}(R\cap F')\Delta([R,R\cap F'])+
\Delta(R')\Delta([R,R\cap F']))\gamma_{5}(R)\\
&\subset&(R\cap F')\cap  (1+\Delta^{4}(R\cap F')
+\Delta(R')\Delta^{2}(R\cap F')
+\Delta(\gamma_{3}(R))\Delta(R\cap F')\\
&&\;\;\;\;\;+\Delta(R')\Delta([R,R\cap F']))
\gamma_{5}(R)[R,R\cap F',R\cap F',R\cap F'].
\end{eqnarray*}
Now since $R\cap F'/R'$ is free-abelian, a use of similar arguments
with left replaced by right and the left $ZR'$-homomorphism
$\phi : Z(R\cap F')\rightarrow ZR'$ implies that
\begin{eqnarray*}
&&(R\cap F')\cap  (1+\Delta^{4}(R\cap F')
+\Delta(R')\Delta^{2}(R\cap F')
+\Delta(\gamma_{3}(R))\Delta(R\cap F')\\
&&\;\;\;\;\;+\Delta(R')\Delta([R,R\cap F']))
\gamma_{5}(R)[R,R\cap F',R\cap F',R\cap F']\\
&=&R'\cap (1+\Delta^{3}(R')+\Delta(R')\Delta([R,R\cap F']))\gamma_{5}(R)\\
&&\;\;\;\;\;[R,R\cap F',R\cap F',R\cap F']\\
&=&\gamma_{5}(R)[R,R\cap F',R\cap F',R\cap F'][\,[R,R\cap F'],[R,R\cap F']\,].
\end{eqnarray*}

\hfill $\Box$

\noindent {\bf \underline{Proof of Theorem A:}} From [6], it follows that 
$$\Delta^3(F)\cap \Delta^2(R)=\Delta^3(R)+\Delta(R)\Delta(R\cap F')+\Delta(R)\Delta([R,R\cap F']),$$
and since $\Delta(R)\mathbb{Z}F$ is a free right $\mathbb{Z}F$-module [3, Proposition I.1.12], we have
$$\Delta^{m}(R)\Delta^{3}(F)\cap \Delta^{m+2}(R)=
\Delta^{m+3}(R)+\Delta^{m+1}(R)\Delta(R\cap F')+\Delta^{m}(R)
\Delta([R,R\cap F'])$$
for all $m\geq 0$. The
natural homomorphism
\[\eta : \Delta^{m+2}(R)\rightarrow
\Delta^{m}(R)\Delta^{2}(F)/\Delta^{m}(R)\Delta^{3}(F)\]
has $ker\,\phi = 
\Delta^{m+3}(R)+\Delta^{m+1}(R)\Delta(R\cap F')+\Delta^{m}(R)
\Delta([R,R\cap F'])$ in view of the above intersection. Thus
$\Delta^{m+2}(R)/
(\Delta^{m+3}(R)+\Delta^{m+1}(R)\Delta(R\cap F')+\Delta^{m}(R)
\Delta([R,R\cap F']))$ is free-abelian. Again, the homomorphism
\[\theta : \gamma_{m+2}(R)\rightarrow
\frac{\Delta^{m+2}(R)}
{\Delta^{m+3}(R)+\Delta^{m+1}(R)\Delta(R\cap F')+\Delta^{m}(R)
\Delta([R,R\cap F'])}\]
defined as $x\rightarrow \overline{(x-1)}\,,x\in \gamma_{m+2}(R)$ has
$ker\,\theta$ equal to
$$\gamma_{m+2}(R)\cap (1+
\Delta^{m+3}(R)+\Delta^{m+1}(R)\Delta(R\cap F')+\Delta^{m}(R)
\Delta([R,R\cap F'])).$$
Therefore
$$\frac{\gamma_{m+2}(R)}{\gamma_{m+2}(R)\cap (1+
\Delta^{m+3}(R)+\Delta^{m+1}(R)\Delta(R\cap F')+\Delta^{m}(R)
\Delta([R,R\cap F']))}$$
is free-abelian for all $m\geq 0$. The proof now follows by putting $m=0,1,2$ in the above group and using Propositions 3.1, 3.2, and 3.3 respectively.


\begin{thebibliography}{[9]}
\bibitem{[1]} Gruenberg K.W., Cohomological Topics in Group Theory, {\em Lecture Notes in Mathematics, Springer, Berlin}, {\bf 143} (1970).
\bibitem{[2]} D.K. Gumber, R. Karan and I. Pal, Some augmentation quotients of integral group rings, {\em Proc. Indian Acad. Sci.} (Math. Sci.) {\bf 118} (2010), 537-546.
\bibitem{[3]} N. Gupta, {\em Free group rings}, Contemporary Math., Amer. Math. Soc. {\bf 66} (1987).
\bibitem{[4]} T. Hurley and S. Sehgal, Groups related to fox
subgroups, {\em Comm. Algebra} {\bf 28} (2000) 1051-1059.
\bibitem{[5]} R. Karan and D. Kumar, Augmentation quotients of free group rings, {\em Algebra Colloq.} {\bf 12} (2005) 597-606.
\bibitem{[6]} R. Karan, D. Kumar and L.R. Vermani, Some intersection theorems and subgroups determined by certain ideals in integral group rings-II, {\em Algebra Colloq.} {\bf 9} (2002), 135-142.
\bibitem{[7]} L.R. Vermani, Augmentation quotients of integral group rings, {\em Groups-Koria'94 (Pusan), de Gruyter, Berlin} (1995) 303-15.
\bibitem{[8]} L.R. Vermani, A. Razdan and R. Karan, Some remarks on subgroups determined by certain ideals in integral group rings, {\em Proc. Indian Acad. Sci. (Math. Sci.)} {\bf 103} (1993), 249-256.   
\bibitem{[9]} K.I. Tahara, L.R. Vermani and Atul Razdan, On generalized third dimension subgroups, {\em Canad. Math. Bull.} {\bf 41} (1998), 109-117.
\end{thebibliography}
\end{document}